\catcode`\@=11

\font\teneufm=eufm10
\font\seveneufm=eufm7
\font\fiveeufm=eufm5
\newfam\eufmfam
\textfont\eufmfam=\teneufm  \scriptfont\eufmfam=\seveneufm
  \scriptscriptfont\eufmfam=\fiveeufm
\def\eufm@{\hexnumber@\eufmfam}

\font\tenmsa=msam10
\font\sevenmsa=msam7
\font\fivemsa=msam5
\font\tenmsb=msbm10
\font\sevenmsb=msbm7
\font\fivemsb=msbm5
\newfam\msafam
\newfam\msbfam
\textfont\msafam=\tenmsa  \scriptfont\msafam=\sevenmsa
  \scriptscriptfont\msafam=\fivemsa
\textfont\msbfam=\tenmsb  \scriptfont\msbfam=\sevenmsb
  \scriptscriptfont\msbfam=\fivemsb

\def\hexnumber@#1{\ifnum#1<10 \number#1\else
 \ifnum#1=10 A\else\ifnum#1=11 B\else\ifnum#1=12 C\else
 \ifnum#1=13 D\else\ifnum#1=14 E\else\ifnum#1=15 F\fi\fi\fi\fi\fi\fi\fi}

\def\msa@{\hexnumber@\msafam}
\def\msb@{\hexnumber@\msbfam}

\mathchardef\gx="2\eufm@78
\mathchardef\gg="2\eufm@67
\mathchardef\gm="2\eufm@6D

\mathchardef\gd="2\eufm@64

\mathchardef\boxdot="2\msa@00
\mathchardef\boxplus="2\msa@01
\mathchardef\boxtimes="2\msa@02
\mathchardef\square="0\msa@03
\mathchardef\blacksquare="0\msa@04
\mathchardef\centerdot="2\msa@05
\mathchardef\lozenge="0\msa@06
\mathchardef\blacklozenge="0\msa@07
\mathchardef\circlearrowright="3\msa@08
\mathchardef\circlearrowleft="3\msa@09
\mathchardef\rightleftharpoons="3\msa@0A
\mathchardef\leftrightharpoons="3\msa@0B
\mathchardef\boxminus="2\msa@0C
\mathchardef\Vdash="3\msa@0D
\mathchardef\Vvdash="3\msa@0E
\mathchardef\vDash="3\msa@0F
\mathchardef\twoheadrightarrow="3\msa@10
\mathchardef\twoheadleftarrow="3\msa@11
\mathchardef\leftleftarrows="3\msa@12
\mathchardef\rightrightarrows="3\msa@13
\mathchardef\upuparrows="3\msa@14
\mathchardef\downdownarrows="3\msa@15
\mathchardef\upharpoonright="3\msa@16

\mathchardef\downharpoonright="3\msa@17
\mathchardef\upharpoonleft="3\msa@18
\mathchardef\downharpoonleft="3\msa@19
\mathchardef\rightarrowtail="3\msa@1A
\mathchardef\leftarrowtail="3\msa@1B
\mathchardef\leftrightarrows="3\msa@1C
\mathchardef\rightleftarrows="3\msa@1D
\mathchardef\Lsh="3\msa@1E
\mathchardef\Rsh="3\msa@1F
\mathchardef\rightsquigarrow="3\msa@20
\mathchardef\leftrightsquigarrow="3\msa@21
\mathchardef\looparrowleft="3\msa@22
\mathchardef\looparrowright="3\msa@23
\mathchardef\circeq="3\msa@24
\mathchardef\succsim="3\msa@25
\mathchardef\gtrsim="3\msa@26
\mathchardef\gtrapprox="3\msa@27
\mathchardef\multimap="3\msa@28
\mathchardef\therefore="3\msa@29
\mathchardef\because="3\msa@2A
\mathchardef\doteqdot="3\msa@2B

\mathchardef\triangleq="3\msa@2C
\mathchardef\precsim="3\msa@2D
\mathchardef\lesssim="3\msa@2E
\mathchardef\lessapprox="3\msa@2F
\mathchardef\eqslantless="3\msa@30
\mathchardef\eqslantgtr="3\msa@31
\mathchardef\curlyeqprec="3\msa@32
\mathchardef\curlyeqsucc="3\msa@33
\mathchardef\preccurlyeq="3\msa@34
\mathchardef\leqq="3\msa@35
\mathchardef\leqslant="3\msa@36
\mathchardef\lessgtr="3\msa@37
\mathchardef\backprime="0\msa@38
\mathchardef\risingdotseq="3\msa@3A
\mathchardef\fallingdotseq="3\msa@3B
\mathchardef\succcurlyeq="3\msa@3C
\mathchardef\geqq="3\msa@3D
\mathchardef\geqslant="3\msa@3E
\mathchardef\gtrless="3\msa@3F
\mathchardef\sqsubset="3\msa@40
\mathchardef\sqsupset="3\msa@41
\mathchardef\trianglerighteq="3\msa@44
\mathchardef\trianglelefteq="3\msa@45
\mathchardef\bigstar="0\msa@46
\mathchardef\between="3\msa@47
\mathchardef\blacktriangledown="0\msa@48
\mathchardef\blacktriangleright="3\msa@49
\mathchardef\blacktriangleleft="3\msa@4A
\mathchardef\blacktriangle="0\msa@4E
\mathchardef\triangledown="0\msa@4F
\mathchardef\eqcirc="3\msa@50
\mathchardef\lesseqgtr="3\msa@51
\mathchardef\gtreqless="3\msa@52
\mathchardef\lesseqqgtr="3\msa@53
\mathchardef\gtreqqless="3\msa@54
\mathchardef\Rrightarrow="3\msa@56
\mathchardef\Lleftarrow="3\msa@57
\mathchardef\veebar="2\msa@59
\mathchardef\barwedge="2\msa@5A
\mathchardef\doublebarwedge="2\msa@5B
\mathchardef\angle="0\msa@5C
\mathchardef\measuredangle="0\msa@5D
\mathchardef\sphericalangle="0\msa@5E
\mathchardef\varpropto="3\msa@5F
\mathchardef\smallsmile="3\msa@60
\mathchardef\smallfrown="3\msa@61
\mathchardef\Subset="3\msa@62
\mathchardef\Supset="3\msa@63
\mathchardef\Cup="2\msa@64

\mathchardef\Cap="2\msa@65

\mathchardef\curlywedge="2\msa@66
\mathchardef\curlyvee="2\msa@67
\mathchardef\leftthreetimes="2\msa@68
\mathchardef\rightthreetimes="2\msa@69
\mathchardef\subseteqq="3\msa@6A
\mathchardef\supseteqq="3\msa@6B
\mathchardef\bumpeq="3\msa@6C
\mathchardef\Bumpeq="3\msa@6D
\mathchardef\lll="3\msa@6E

\mathchardef\ggg="3\msa@6F

\mathchardef\circledS="0\msa@73
\mathchardef\pitchfork="3\msa@74
\mathchardef\dotplus="2\msa@75
\mathchardef\backsim="3\msa@76
\mathchardef\backsimeq="3\msa@77
\mathchardef\complement="0\msa@7B
\mathchardef\intercal="2\msa@7C
\mathchardef\circledcirc="2\msa@7D
\mathchardef\circledast="2\msa@7E
\mathchardef\circleddash="2\msa@7F
\def\ulcorner{\delimiter"4\msa@70\msa@70 }
\def\urcorner{\delimiter"5\msa@71\msa@71 }
\def\llcorner{\delimiter"4\msa@78\msa@78 }
\def\lrcorner{\delimiter"5\msa@79\msa@79 }
\def\yen{\mathhexbox\msa@55 }
\def\checkmark{\mathhexbox\msa@58 }
\def\circledR{\mathhexbox\msa@72 }
\def\maltese{\mathhexbox\msa@7A }
\mathchardef\lvertneqq="3\msb@00
\mathchardef\gvertneqq="3\msb@01
\mathchardef\nleq="3\msb@02
\mathchardef\ngeq="3\msb@03
\mathchardef\nless="3\msb@04
\mathchardef\ngtr="3\msb@05
\mathchardef\nprec="3\msb@06
\mathchardef\nsucc="3\msb@07
\mathchardef\lneqq="3\msb@08
\mathchardef\gneqq="3\msb@09
\mathchardef\nleqslant="3\msb@0A
\mathchardef\ngeqslant="3\msb@0B
\mathchardef\lneq="3\msb@0C
\mathchardef\gneq="3\msb@0D
\mathchardef\npreceq="3\msb@0E
\mathchardef\nsucceq="3\msb@0F
\mathchardef\precnsim="3\msb@10
\mathchardef\succnsim="3\msb@11
\mathchardef\lnsim="3\msb@12
\mathchardef\gnsim="3\msb@13
\mathchardef\nleqq="3\msb@14
\mathchardef\ngeqq="3\msb@15
\mathchardef\precneqq="3\msb@16
\mathchardef\succneqq="3\msb@17
\mathchardef\precnapprox="3\msb@18
\mathchardef\succnapprox="3\msb@19
\mathchardef\lnapprox="3\msb@1A
\mathchardef\gnapprox="3\msb@1B
\mathchardef\nsim="3\msb@1C
\mathchardef\napprox="3\msb@1D
\mathchardef\nsubseteqq="3\msb@22
\mathchardef\nsupseteqq="3\msb@23
\mathchardef\subsetneqq="3\msb@24
\mathchardef\supsetneqq="3\msb@25
\mathchardef\subsetneq="3\msb@28
\mathchardef\supsetneq="3\msb@29
\mathchardef\nsubseteq="3\msb@2A
\mathchardef\nsupseteq="3\msb@2B
\mathchardef\nparallel="3\msb@2C
\mathchardef\nmid="3\msb@2D
\mathchardef\nshortmid="3\msb@2E
\mathchardef\nshortparallel="3\msb@2F
\mathchardef\nvdash="3\msb@30
\mathchardef\nVdash="3\msb@31
\mathchardef\nvDash="3\msb@32
\mathchardef\nVDash="3\msb@33
\mathchardef\ntrianglerighteq="3\msb@34
\mathchardef\ntrianglelefteq="3\msb@35
\mathchardef\ntriangleleft="3\msb@36
\mathchardef\ntriangleright="3\msb@37
\mathchardef\nleftarrow="3\msb@38
\mathchardef\nrightarrow="3\msb@39
\mathchardef\nLeftarrow="3\msb@3A
\mathchardef\nRightarrow="3\msb@3B
\mathchardef\nLeftrightarrow="3\msb@3C
\mathchardef\nleftrightarrow="3\msb@3D
\mathchardef\divideontimes="2\msb@3E
\mathchardef\varnothing="0\msb@3F
\mathchardef\nexists="0\msb@40
\mathchardef\mho="0\msb@66
\mathchardef\thorn="0\msb@67
\mathchardef\beth="0\msb@69
\mathchardef\gimel="0\msb@6A
\mathchardef\daleth="0\msb@6B
\mathchardef\lessdot="3\msb@6C
\mathchardef\gtrdot="3\msb@6D
\mathchardef\ltimes="2\msb@6E
\mathchardef\rtimes="2\msb@6F
\mathchardef\shortmid="3\msb@70
\mathchardef\shortparallel="3\msb@71
\mathchardef\smallsetminus="2\msb@72
\mathchardef\thicksim="3\msb@73
\mathchardef\thickapprox="3\msb@74
\mathchardef\approxeq="3\msb@75
\mathchardef\succapprox="3\msb@76
\mathchardef\precapprox="3\msb@77
\mathchardef\curvearrowleft="3\msb@78
\mathchardef\curvearrowright="3\msb@79
\mathchardef\digamma="0\msb@7A
\mathchardef\varkappa="0\msb@7B
\mathchardef\hslash="0\msb@7D
\mathchardef\hbar="0\msb@7E
\mathchardef\backepsilon="3\msb@7F
\def\Bbb{\ifmmode\let\next\Bbb@\else
 \def\next{\errmessage{Use \string\Bbb\space only in math mode}}\fi\next}
\def\Bbb@#1{{\Bbb@@{#1}}}
\def\Bbb@@#1{\fam\msbfam#1}

\catcode`\@=\active

\def\eps{{\epsilon}}

\def\<{\langle}
\def\>{\rangle}

\def\tens{\mathop{\otimes}}
\def\la{{\triangleright}}\def\ra{{\triangleleft}}

\def\proof{\goodbreak\noindent{\bf Proof\quad}}

\def\text#1{{\rm #1}}
\def\note#1{}

\documentclass[12pt]{article}
\newtheorem{lemma}{Lemma}[section]
\newtheorem{propos}[lemma]{Proposition}
\newtheorem{example}[lemma]{Example}

\newtheorem{defin}[lemma]{Definition}

\begin{document}

\newpage

\hfill Swan-Maths-01/1

\begin{center} {\Large Two-forms and Noncommutative Hamiltonian dynamics}
\\ \baselineskip 13pt{\ }
{\ }\\  E. J. Beggs\\
{\ } \\ Department of Mathematics\\
University of Wales Swansea\\ Wales SA2 8PP
\end{center}

\begin{quote}
\noindent{\bf Abstract.} 
In this paper we extend the standard differential geometric theory of Hamiltonian dynamics to noncommutative
spaces, beginning with symplectic forms. Derivations on the algebra are used instead of vector fields, and
interior products and Lie derivatives with respect to derivations are discussed. Then the Poisson bracket of
certain algebra elements can be defined by a choice of closed 2-form. Examples are given using the noncommutative
torus, the Cuntz algebra, the algebra of matrices, and the algebra of matrix valued functions on $\Bbb R^2$. 

\end{quote}

\section{Introduction}
We begin from the usual definition of a differential calculus on a noncommutative algebra. The 
derivations on the algebra are used to substitute for the vector fields in the commutative case. 
Then we can define an interior product and Lie derivative, and prove noncommutative analogues 
of the standard results in classical differential geometry. The proofs are almost identical to the
classical ones. The caveat is that we must consider only those derivations which are compatible with the relations
in the differential structure. 

From a closed 2-form $\omega$ we define a map
$\tilde\omega$ from the collection of derivations to the 1-forms. Classically this would be a 1-1
correspondence if $\omega$ was nondegenerate. We shall just assume that $\tilde\omega$ is 1-1, and as a result 
only certain elements of the algebra (called Hamiltonian elements) will correspond to derivations. 
The variety of examples of differential calculi on noncommutative algebras means that to insist on comparable
sizes for the set of derivations and the 1-forms would be over restrictive, and the
fact that $\tilde\omega$ might not be onto will not cause major problems. For Hamiltonian elements
we can define a Poisson bracket, which is antisymmetric and satisfies the Jacobi identity.
It may be suprising that a noncommutative differential geometry can have antisymmetric Poisson brackets. 
However the reader should note that we do not achieve this by imposing any sort of asymmetry on the
differential forms, but by imposing asymmetry on the interior product $\lrcorner$, by insisting that
it is a signed derivation. 

Then there are the examples. For the noncommutative torus $\Bbb T_\rho^2$ (with $uv=e^{2\pi i\rho}\,vu$)
we take $\omega=u^{-1}\,
du\,dv\, v^{-1}$.  If $\rho$ is rational there is a class of Hamiltonian elements with
 non-zero Poisson brackets, although
the Hamiltonian elements commute in the algebra multiplication. For the matrix algebra $M_n(\Bbb R)$
with $\omega=\sum_{ij} dE_{ij}\, dE_{ij}$ we see that the antisymmetric matrices are Hamiltonian elements,
and the corresponding derivations are just the adjoint maps of the antisymmetric matrices. The Poisson
bracket is just the matrix commutator. The Cuntz algebra
${\cal O}_n$ provides another example. 
We conclude by examining Hamiltonian dynamics on the algebra of matrix valued functions on $\Bbb R^2$.

It is time for a public health warning: Differential calculi on $C^*$ algebras frequently require passing to a
smaller `smooth' subalgebra to make things work. For an example see the 
 calculation of cyclic cohomology in \cite{connesIHES}. We shall work purely algebraically
in what follows, without worrying
about the topology.

The author would like to thank Tomasz Brzezi\'nski for much useful advice.

\section{Differential calculus and derivations}

\begin{defin} A differential calculus on an algebra $A$ is a collection of $A$-bimodules
$\Omega^n$ for $n\ge 0$ and a signed derivation $d:\Omega^n\to\Omega^{n+1}$,
 i.e. $d(\omega\tau)=d(\omega)\,\tau+(-1)^{|\omega|}\,\omega\, d(\tau)$. 
Here $|\omega|=n$ if $\omega\in\Omega^n$. We set $\Omega^0=A$
and suppose that the subspace spanned by elements of the form $\omega\, db$ ( for all
$\omega\in\Omega^n$ and $b\in A$)
is dense in $\Omega^{n+1}$. We also impose $d^2=0$. 
\end{defin}

\begin{defin} 
Define $V$ to be a vector space of derivations on the algebra $A$, i.e.\ 
$\theta(ab)=\theta(a)b+a\theta(b)$ for  $\theta\in V$. We suppose that $V$ is closed 
under the commutator $[\theta,\phi]=\theta\phi-\phi\theta$. 
 This will take the place of the
vector fields in commutative differential geometry.
\end{defin} 

\smallskip Now that we have the analogue of vector fields, we can define the following operations. 

\begin{defin} 
We define the evaluation map or `interior product' $\lrcorner:V\tens\Omega^1 \to A$
 by $\theta\lrcorner da=\theta(a)$. 
 To be consistent with the rule $d(ab)=da\, b+a\, db$ we set
$\theta\lrcorner(da\, b)=\theta(a)\, b$ and $\theta\lrcorner(a\,db)=a\,\theta(b)$. Extend
this definition recursively to $\lrcorner:V\tens \Omega^{n+1}\to\Omega^n$ as a signed derivation, i.e.\ 
$\theta\lrcorner(\omega\tau)=(\theta\lrcorner \,\omega)\,\tau\,+\, (-1)^{|\omega|}\,\omega\, (\theta\lrcorner\tau)$. 
\end{defin} 

\begin{defin} 
We define the Lie derivative ${\cal L}_\theta:\Omega^1\to\Omega^1$ in the direction $\theta\in V$
 of a 1-form by ${\cal L}_\theta(da)=d(\theta(a))$, and extend it as a derivation, i.e.\
$\theta(a\,db)=\theta(a)\, db\,+\, a\, d(\theta(b))$.
This is compatible with the rule $d(ab)=da\, b+a\, db$.  Now we extend the definition to
${\cal L}_\theta:\Omega^n\to\Omega^n$ as a derivation, i.e.\ ${\cal L}_\theta(\omega\tau)
=\omega{\cal L}_\theta(\tau)+{\cal L}_\theta(\omega)\tau$. 
\end{defin} 

\smallskip The problem with these operations is that they may not be well defined, that is there may be a linear
combination of elements of the form $a\,db\in \Omega^1$ which vanishes, but for which the corresponding sum of
$\theta
\lrcorner (a\,db)$ or ${\cal L}_\theta(a\,db)$ would not be zero. If the differential calculus is given in terms
of generators and relations, we must check that the interior product and the Lie derivative vanish on all the
relations. If necessary we must restrict the set $V$ of derivations so that these operations are well defined.
  In what follows, we assume that these operations are well defined.

\begin{propos}\label{propa} For all $\theta\in V$ and $\omega\in \Omega^n$,
$
d(\theta\lrcorner \,\omega)\,+\, \theta\lrcorner(d\omega)\, =\, {\cal L}_\theta(\omega)
$.
\end{propos}
\proof By induction on the degree of $\omega$. The statement is true for 0-forms (elements of $A$).
Now we suppose that the statement is true for $n$-forms, and take $\omega\in\Omega^n$
and $b\in A$. 
\begin{eqnarray*} 
d(\theta\lrcorner(\omega\, db))&=& d(\, (\theta\lrcorner \,\omega)\, db\,+\, (-1)^n\,\omega\, \theta(b)\,) \cr
&=& d(\theta\lrcorner \,\omega)\, db\,+\, (-1)^n\,d\omega\, \theta(b)\,+\,\omega\,d\theta(b)\ ,\cr
\theta\lrcorner(d(\omega\, db))&=& \theta\lrcorner(d\omega\, db)\ =\ 
(\theta\lrcorner d\omega)\, db\,+\, (-1)^{n+1}\, d\omega\,\theta(b)\ .
\end{eqnarray*} 
Now add these together to get
\begin{eqnarray*} 
d(\theta\lrcorner(\omega\, db))\,+\, \theta\lrcorner(d(\omega\, db))&=& 
(d(\theta\lrcorner \,\omega)\,+\, \theta\lrcorner d\omega)\, db\,+\, \omega\,d\theta(b) \cr
&=& {\cal L}_\theta(\omega)\, db\,+\,  \omega\,d\theta(b)\ =\ {\cal L}_\theta(\omega\, db)\ .\square
\end{eqnarray*}

\begin{propos} For all $\theta\in V$ and $\omega\in \Omega^n$,
$
d{\cal L}_\theta(\omega)\, =\, {\cal L}_\theta(d\omega)
$.
\end{propos}
\proof By induction on the degree of $\omega$. The statement is true for 0-forms (elements of $A$).
Now we suppose that the statement is true for $n$-forms, and take $\omega\in\Omega^n$
and $b\in A$. 
\begin{eqnarray*} 
d{\cal L}_\theta(\omega\, db) &=& d(\ ({\cal L}_\theta(\omega))\, db\,+\, \omega\, d\theta(b)\ ) \cr
 &=& (d{\cal L}_\theta(\omega))\, db\,+\, d\omega\, d\theta(b) \cr 
&=& {\cal L}_\theta(d\omega)\, db\,+\, d\omega\, {\cal L}_\theta(db) \ =\ 
{\cal L}_\theta(d\omega\, db)\ =\ {\cal L}_\theta(d(\omega\, db))\ .\square
\end{eqnarray*} 

\begin{propos}\label{propc} For all $\theta,\phi\in V$ and $\omega\in \Omega^n$,
$
{\cal L}_\phi(\theta\lrcorner \,\omega)\, =\, \theta\lrcorner {\cal L}_\phi(\omega)\,+\,
[\phi,\theta]\lrcorner \,\omega
$.
\end{propos}
\proof By induction on the degree of $\omega$. The statement is true for 0-forms (elements of $A$).
Now we suppose that the statement is true for $n$-forms, and take $\omega\in\Omega^n$
and $b\in A$. 
Then 
\begin{eqnarray*} 
{\cal L}_\phi(\theta\lrcorner (\omega\,db)) &=& {\cal L}_\phi((\theta\lrcorner \,\omega)\,db\,+\,
(-1)^n\,\omega\, \theta(b)) \cr
&=& 
{\cal L}_\phi(\theta\lrcorner \,\omega)\,db\,+\, (\theta\lrcorner \,\omega)\, d\phi(b)
\,+\,
(-1)^n\,{\cal L}_\phi(\omega)\, \theta(b)\,+\,
(-1)^n\,\omega\, \phi\theta(b)   \ ,\cr
\theta\lrcorner({\cal L}_\phi(\omega\,db)) &=& \theta\lrcorner({\cal L}_\phi(\omega)\,db\,+\,
\omega\,d\phi(b)) \cr
&=& 
(\theta\lrcorner{\cal L}_\phi(\omega))\,db\,+\,(-1)^n\,{\cal L}_\phi(\omega)\,\theta(b)
\,+\, (\theta\lrcorner \,\omega)\,d\phi(b)\,+\,(-1)^n\,
\omega\,\theta\phi(b)\ ,
\end{eqnarray*} 
and on subtraction we get
\begin{eqnarray*} 
{\cal L}_\phi(\theta\lrcorner \,\omega\,db)\,-\, \theta\lrcorner({\cal L}_\phi(\omega\,db)) &=&
{\cal L}_\phi(\theta\lrcorner \,\omega)\,db\,-\,
(\theta\lrcorner{\cal L}_\phi(\omega))\,db\,+\, (-1)^n\,\omega\, [\phi,\theta](b) \cr
 &=&
([\phi,\theta]\lrcorner \,\omega)\,db\,+\, (-1)^n\,\omega\, [\phi,\theta](b)\cr
&=& [\phi,\theta]\lrcorner (\omega\,db)\ .\qquad\square
\end{eqnarray*} 

\begin{propos}\label{propb} For all $\theta,\phi\in V$ and $\omega\in \Omega^n$,
\[
\phi\lrcorner(\theta\lrcorner \,\omega)\ =\ -\, \theta\lrcorner(\phi\lrcorner \,\omega)\ .
\]
\end{propos}
\proof By induction on the degree of $\omega$. The statement is true for 0-forms (elements of $A$).
Now we suppose that the statement is true for $n$-forms, and take $\omega\in\Omega^n$
and $b\in A$. 
\begin{eqnarray*} 
\phi\lrcorner(\theta\lrcorner(\omega\, db)) &=& \phi\lrcorner((\theta\lrcorner \,\omega)\, db\,+\,(-1)^n\,\omega\,\theta(b))
\cr
&=& (\phi\lrcorner(\theta\lrcorner \,\omega))\, db\,+\, (-1)^{n-1}\, (\theta\lrcorner \,\omega)\, \phi(b)
\,+\, (-1)^n\,(\phi\lrcorner \,\omega)\,\theta(b)\ .
\end{eqnarray*} 
Now just add this formula to the one with $\phi$ and $\theta$ swapped to get zero.\quad$\square$

\begin{propos} \label{proppp} For all $\theta,\phi\in V$ and $\omega\in \Omega^n$,
\[
{\cal L}_\theta {\cal L}_\phi (\omega)\,-\,  {\cal L}_\phi {\cal L}_\theta (\omega)
\ =\ {\cal L}_{[\theta,\phi]}(\omega)\ .
\]
\end{propos}
\proof By induction on the degree of $\omega$. The statement is true for 0-forms (elements of $A$).
Now we suppose that the statement is true for $n$-forms, and take $\omega\in\Omega^n$
and $b\in A$.
\begin{eqnarray*} 
{\cal L}_\theta {\cal L}_\phi (\omega\,db)&=& {\cal L}_\theta ({\cal L}_\phi (\omega)\, db\,+\,
\omega\,d\phi(b)) \cr
&=& {\cal L}_\theta {\cal L}_\phi (\omega)\, db\,+\,{\cal L}_\theta(\omega)\,d\phi(b)
\,+\,{\cal L}_\phi(\omega)\,d\theta(b) \,+\,  \omega\,d\theta\phi(b)\ ,
\end{eqnarray*} 
and if we swap $\theta$ and $\phi$ and subtract we get
\[
({\cal L}_\theta {\cal L}_\phi\,-\, {\cal L}_\phi {\cal L}_\theta) (\omega\,db)\ =\ 
({\cal L}_\theta {\cal L}_\phi\,-\, {\cal L}_\phi {\cal L}_\theta) (\omega)\,db\,+\,
\omega\, d[\theta,\phi](b)\ =\ {\cal L}_{[\theta,\phi]} (\omega\,db)\ .\square
\]

\section{Hamiltonian dynamics}
In this section we take a specified $\omega\in \Omega^2$ which is closed, i.e.\ $d\omega=0$. 
From proposition \ref{propa} we see that ${\cal L}_\theta(\omega)=0$ if and only
if $d(\theta\lrcorner \,\omega)=0$.  Take the subset $V^\omega$ to consist of those
$\theta\in V$ for which ${\cal L}_\theta(\omega)=0$, and $Z^1$ to be the set of closed 1-forms. 
Define the map $\tilde\omega:V^\omega\to Z^1$ by
$\tilde\omega(\theta)=\theta\lrcorner \,\omega$.  We say that $\omega$ is nonsingular if $\tilde\omega$ 
is 1-1, and we suppose this for the rest 
of the section.

\begin{defin} We say that $a\in A$ is a Hamiltonian element if $da\in Z^1$ is in the 
image of $\tilde\omega$. If $a$ is Hamiltonian, we define $X_a\in V^\omega$ by 
$\tilde\omega(X_a)=da$. If both $a$ and $b$ are Hamiltonian, we define
their Poisson bracket by $\{a,b\}=X_a\lrcorner(db)=
X_a(b)\in A$. 
\end{defin}

\begin{propos} If  both $a$ and $b$ are Hamiltonian, then
$\{a,b\}=-\{b,a\}$, i.e.\ the Poisson bracket is antisymmetric.
\end{propos} 
\proof From proposition \ref{propb},
\[
X_a\lrcorner(X_b \lrcorner \,\omega)\ =\ X_a\lrcorner db\ =\ X_a(b)\ =\ -\, 
X_b\lrcorner(X_a \lrcorner \,\omega)\ =\ -\, X_b(a)\ .\quad\square
\]

\begin{propos}\label{propd} If  both $a$ and $b$ are Hamiltonian, then
$\{a,b\}$ is Hamiltonian, and further $X_{\{a,b\}}=[X_a,X_b]$. 
\end{propos} 
\proof First $[X_a,X_b]\in V$ as $V$ is closed under commutator. Then ${\cal L}_{[X_a,X_b]}=0$ by
\ref{proppp}, so $[X_a,X_b]\in V^\omega$. Finally from proposition \ref{propc},
\[ 
[X_a,X_b] \lrcorner \,\omega \ =\ {\cal L}_{X_a}(X_b \lrcorner \,\omega)\,-\, X_b \lrcorner
{\cal L}_{X_a}(\omega)\ =\ {\cal L}_{X_a}(db)\ =\ d X_a(b)\ =\ d\{a,b\}\ .\quad\square
\]

\begin{propos} If $a,b$ and $c$ are Hamiltonian, then
$\{c,\{a,b\}\}\,+\,\{b,\{c,a\}\}\,+\,\{a,\{b,c\}\}\,=\,0$, i.e.\ the 
Poisson bracket satisfies the Jacobi identity.
\end{propos} 
\proof By using proposition \ref{propc},
\begin{eqnarray*} 
X_c\{a,b\} &=& X_c \lrcorner d\{a,b\}\ =\  X_c \lrcorner d X_{a}(b)\ =\ X_c \lrcorner {\cal L}_{X_a}(db)
   \cr
 &=& 
{\cal L}_{X_a}(X_c  \lrcorner db)\,-\,
[X_a,X_c] \lrcorner db\ .
\end{eqnarray*} 
From this we deduce, using \ref{propd},
\begin{eqnarray*} 
\{c,\{a,b\}\}\,+\, \{\{a,c\},b\}&=&
 {\cal L}_{X_a}(X_c  \lrcorner db)\ =\ {\cal L}_{X_a}\{c,b\}\ =\ \{a,\{c,b\}\}\ .\square
\end{eqnarray*}

\begin{propos} If $a,b,c$ and $bc$ are Hamiltonian, then
$\{a,bc\}\,=\,\{a,b\}\,c\,+\,b\,\{a,c\}$, i.e.\ the 
Poisson bracket is a derivation.
\end{propos} 
\proof Use the result $\{a,bc\}\,=\,X_a(bc)$, where $X_a$ is a derivation.\quad$\square$

\smallskip

Now we can formally extend a derivation to an automorphism by the following procedure
(we make no attempt to verify convergence): If $\theta$ is a derivation
on the algebra $A$, there is an action of $(\Bbb R,+)$ by automorphisms on $A$ by
$a\mapsto \exp(t\theta)a=a(t)$ for $a\in A$ and $t\in \Bbb R$. Then we get the usual relation
for the time derivatives of functions $a(t)\in A$ generated by a Hamiltonian
$b\in A$ and the Poisson bracket:
\[
\dot a(t)\ =\ X_b(a(t))\ =\ \{b,a(t)\}\ .
\]
(Note that strictly we should stop at $X_b(a(t))$ in the case where $a$ is not Hamiltonian, as we did
not define the Poisson brackets for non-Hamiltonian elements.)

\section{Example: the noncommutative torus}
We take the algebra $\Bbb T^2_\rho$ generated by invertible elements $u$ and $v$, subject to the conditions
$uv=qvu$, where $q=e^{2\pi \rho}$ is a unit norm complex number. This can be completed to form a $C^*$ algebra,
or a  smooth algebra \cite{connesIHES,eeirr}, but we will not consider such completions here. The simplest
differential  calculus on $\Bbb T^2_q$ \cite{bdr} is generated by $\{u,v,du,dv\}$, subject to the relations
\begin{eqnarray} \label{nctrel}
du\,dv=-q\, dv\,du\ &,& u\, dv=q\, dv\, u\ ,\quad v\, du=q^{-1}du\, v\ ,\cr
[u,du]=[v,dv]=0\ &,& (du)^2=(dv)^2=0\ .
\end{eqnarray} 
As there are no non-zero 3-forms, all 2-forms are closed. 

We will now try to carry out the construction given in the previous sections. 
Note that a derivation $\theta$ is
uniquely specified by giving $\theta(u)$ and $\theta(v)$. This is because we can deduce
$\theta(u^{-1})=-\,u^{-1}\,\theta(u)\,u^{-1}$ from the relation $\theta(u\,u^{-1})=\theta(1)=0$,
and likewise for $v^{-1}$. 
The set of derivations $V$ which we use must be consistent with the relations on the
algebra and on the differential structure.

\begin{propos} Suppose that we are
in the rational case where $q$ has order $p$.
Then the derivations $\theta$ consistent with the differential structure
(\ref{nctrel}) are of the form
\[
\theta(u)\ =\ \sum_{s,t\in\Bbb Z} b_{ts}\, u^{1+sp}\, v^{tp}\ ,\quad
\theta(v)\ =\ \sum_{s,t\in\Bbb Z} c_{ts}\, u^{sp}\, v^{1+tp}\ ,
\]
for some constants $b_{ts}$, $c_{ts}\in\Bbb C$
\end{propos} 
\proof By applying $\theta\lrcorner $ to $[u,du]=[v,dv]=0$ we see that $[u,\theta(u)]=[v,\theta(v)]=0$. 
By applying $\theta\lrcorner $ to $u\, dv=q\, dv\, u$ we see that $u\, \theta(v)=q\, \theta(v)\, u$. 
By applying $\theta\lrcorner $ to $v\, du=q^{-1}du\, v$ we see that $v\, \theta(u)=q^{-1}\theta(u)\, v$. 
By combining these we see that the consistent derivations are those for which
are of the form above. 
Now we check with the algebra relation $\theta(uv)=q\theta(vu)$ to get
\[
\sum_{s,t\in\Bbb Z} c_{ts}\, u^{1+sp}\, v^{1+tp}\ +\ 
\sum_{s,t\in\Bbb Z} b_{ts}\, u^{1+sp}\, v^{1+tp}\ =\ 
q\, \sum_{s,t\in\Bbb Z} c_{ts}\, u^{sp}\, v^{1+tp}\, u\ +\ 
q\, v\, \sum_{s,t\in\Bbb Z} b_{ts}\, u^{1+sp}\, v^{tp}\ ,
\]
which is automatically satisfied.\quad$\square$

\begin{example} 
Set $\omega=u^{-1}\,du\, dv\,v^{-1}$,
and suppose that we are
in the rational case where $q$ has order $p$. Then for $\theta\in V$ we have $\theta\lrcorner \,\omega=u^{-1}\,\theta(u)\,
dv\,
v^{-1}\,-\, u^{-1}\,du\,\theta(v)\,v^{-1}$, so if we know $\theta\lrcorner \,\omega$ we can recover $\theta(u)$ and
$\theta(v)$ uniquely, so $\tilde\omega$ is 1-1. 
Proceeding on the assumption that $a\in \Bbb T^2_\rho$ is a Hamiltonian element, we set
\[
a\ =\ \sum_{nm} a_{nm}\, u^n\, v^m
\]
for some numbers $a_{nm}$. We now examine the equation 
\begin{eqnarray*} 
X_a\lrcorner(u^{-1}\,du\, dv\,v^{-1}) &=&  da \ =\ 
 \sum_{nm} \Big(  n\,a_{nm}\, u^{n-1}\, du\, v^m\,+\,m\, a_{nm}\, u^n\, v^{m-1}\, dv\Big)\ ,
\end{eqnarray*} 
and deduce that
\[
X_a(u)\ =\ \sum_{nm} m\,  a_{nm}\, u^{n+1}\, v^{m} \quad{\rm and}\quad
X_a(v)\ =\ -\sum_{nm} n\,a_{nm}\, u^{n}\, v^{m+1}\ .
\]
For $X_a$ to be a derivation consistent with our given differential structure,
 we can only have nonzero $a_{nm}$
when $n$ and $m$ are multiples of $p$. The Hamiltonian functions are linear combinations
of elements of the form $u^{sp}\, v^{tp}$ for $t,s\in\Bbb Z$. 
The corresponding derivations are
\[
X_{u^{sp}\, v^{tp}}(u)\ =\ tp\, u^{sp+1}\, v^{tp} \quad{\rm and}\quad
X_{u^{sp}\, v^{tp}}(v)\ =\ -\,sp\,u^{sp}\, v^{tp+1}\ ,
\]
and the Poisson brackets are given by
\[
\{u^{sp}\, v^{tp}\,,\,u^{s'p}\, v^{t'p}\}\ =\ 
(t\,s'\,-\, t'\,s)\,p^2\, u^{(s+s')p}\, v^{(t+t')p}\ .
\]
Note that the Hamiltonian elements in this case are
 exactly the central elements of the algebra. 

\end{example}

\section{Example: the algebra of matrices}
We take $A=M_n(\Bbb R)$, and then we define $\Omega^1$ to be the kernel of the multiplication map 
$\mu:M_n\tens M_n\to M_n$ \cite{kar,connesIHES}. The map $d:\Omega^0=A\to \Omega^1$ is defined as
$da=1\tens a-a\tens 1$. Also
\[
\Omega^2\ =\ \Big\{ \tau\in M_n\tens M_n\tens M_n:(\mu\tens {\rm id})(\tau)=
({\rm id}\tens \mu)(\tau)=0\Big\}\ .
\]
The map $d:\Omega^1\to\Omega^2$ is defined by
$d(a\tens b)=1\tens a\tens b-a\tens 1\tens b+a\tens b\tens 1$. In case the reader is 
concerned that the interior product is not well defined on this model of the differential calculus,
note that $\theta\lrcorner:\Omega^1\to \Omega^0=M_n$ is 
$\theta\lrcorner(a\tens b)=a\,\theta(b)$, and $\theta\lrcorner:\Omega^2\to \Omega^1$ is 
$\theta\lrcorner(a\tens b\tens c)=a\,\theta(b)\tens c\,-\,a\tens b\,\theta(c)$. 
The Lie derivative ${\cal L}_\theta:\Omega^1\to\Omega^1$ is given by
${\cal L}_\theta(a\tens b)\,=\, \theta(a)\tens b\,+\,a\tens\theta(b)$, and 
${\cal L}_\theta:\Omega^2\to\Omega^2$ is given by
${\cal L}_\theta(a\tens b\tens c)\,=\, \theta(a)\tens b\tens c\,+\,a\tens\theta(b)\tens c
\,+\,a\tens b\tens \theta(c)$.

Set $\omega=\frac12\sum_{ij}dE_{ij}\, dE_{ij}$, where $E_{ij}$ is the matrix with 1 in the row
$i$ column $j$ position and 0 elsewhere. Take a derivation $\theta$ on 
$M_n(\Bbb R)$ given by coefficients $\Theta_{klij}\in\Bbb R$:
\[
\theta(E_{ij})\ =\ \sum_{kl}\Theta_{klij}\, E_{kl}\ .
\]
Then we calculate
\begin{eqnarray*} 
\theta\lrcorner \,\omega &=& \frac12\sum_{ijkl}\Theta_{klij}\,\Big(E_{kl}\,dE_{ij}\,-\, dE_{ij}\, E_{kl}\Big)
\cr
 &=& \frac12\sum_{ijkl}\Theta_{klij}\,\Big(E_{kl}\tens E_{ij}\,-\, E_{kl}\, E_{ij}\tens 1\,-\,
1\tens E_{ij}\, E_{kl}\,+\, E_{ij}\tens E_{kl}
\Big)
\cr
 &=& \frac12\sum_{ijkl}\Big(\Theta_{klij}+\Theta_{ijkl}\Big)\,E_{kl}\tens E_{ij}\,-\,
\frac12\sum_{ijkl} \Theta_{klij}\, E_{kl}\, E_{ij}\tens 1\,-\,
\frac12\sum_{ijkl} \Theta_{ijkl}\,1\tens E_{kl}\, E_{ij}\ .
\end{eqnarray*} 
Given a matrix $S\in M_n(\Bbb R)$, take the adjoint map ${\rm ad}_S(C)=[S,C]$, which is a derivation. 
Now
\[
{\rm ad}_S(E_{ij})\ =\ [S,E_{ij}]\ =\ \sum_k S_{ki}\,E_{kj}\,-\, \sum_l E_{il}\, S_{jl}\ ,
\]
so the coefficients corresponding to $\theta={\rm ad}_S$ are $\Theta_{klij}=
 S_{ki}\delta_{jl}\,-\,S_{jl}\delta_{ik}$. If £S£ is an antisymmetric matrix, then
$\Theta_{klij}\,+\,\Theta_{ijkl}\,=\,0$, so 
\begin{eqnarray*} 
{\rm ad}_S\lrcorner \,\omega &=& 
\frac12\sum_{ijkl} \Theta_{klij}\, \Big(1\tens E_{kl}\, E_{ij}\,-\,E_{kl}\, E_{ij}\tens 1\Big)
\cr &=& 
\frac12\sum_{ijkl}\Big(S_{ki}\delta_{jl}\,-\,S_{jl}\delta_{ik}\Big)\, \delta_{li}\,
\Big(1\tens E_{kj}\,-\,E_{kj}\tens 1\Big)\ =\ 1\tens S\,-\, S\tens 1\ .
\end{eqnarray*} 
Now we see that ${\rm ad}_S\lrcorner \,\omega=dS$, so the antisymmetric matrices are Hamiltonian,
and $X_S={\rm ad}_S$. If $S$ and $T$ are antisymmetric, then
$\{S,T\}={\rm ad}_S(T)=[S,T]$, so the Poisson bracket is just the commutator.

\section{Example: the Cuntz algebra} 
The Cuntz algebra ${\cal O}_n$ \cite{cu1} is the unital $C^*$ algebra with
 $n$ generators $\{s_1,\dots,s_n\}$ and relations
\[
s_i^*\, s_j\ =\ \delta_{ij}\ ,\quad 
\sum_{i=1,\dots,n} s_i\, s_i^*\ =\ 1\ .
\]
Linear combinations of the form $s_\mu\,s_\nu^*$ are dense in the algebra,
where $\mu$ and $\nu$ are words for the alphabet $\{1,\dots,n\}$. 
For example if $\mu=12$ and $\nu=123$, then
$s_\mu\,s_\nu^*=s_1s_2s_3^*s_2^*s_1^*$. 

Now we have to decide what differential calculus to equip ${\cal O}_n$ with. The 
forms would be generated by $s_i$, $s_i^*$, $ds_i$ and $ds_i^*$. We must have relations
given by applying $d$ to the relations for the algebra, i.e.\
\begin{eqnarray} \label{cudifrel}
ds_i^*\, s_j\,+\, s_i^*\, ds_j\ =\ 0\ ,\quad 
\sum_{i=1,\dots,n} (ds_i\, s_i^*\,+\, s_i\, ds_i^*)\ =\ 0\ .
\end{eqnarray} 

If $u$ is any unitary in ${\cal O}_n$, then the map
$s_i\mapsto u\,s_i$ extends to a unital *-endomorphism $\alpha_u$ of ${\cal O}_n$. 
Conversely, suppose that $\alpha$ is a
unital *-endomorphism of ${\cal O}_n$. 
If we define $u=\sum \alpha(s_i)\, s_i^*$, we see that 
$u$ is a unitary in ${\cal O}_n$, and that $\alpha(s_i)=u\,s_i$. 
We shall use this to define a derivation on ${\cal O}_n$ by taking the infinitesimal
version of this construction. 
For $h\in {\cal O}_n$ we define a derivation by $\theta_h(s_i)=h\, s_i$ and 
$\theta_h(s_i^*)=-s_i^*\, h$. If this were to be a *-derivation, we would find
that $h$ had to be antiHermitian, but we shall not suppose this. Now we should check
 that these derivations
preserve the relations (\ref{cudifrel}):
\begin{eqnarray*}
\theta_h\lrcorner (ds_i^*\, s_j\,+\, s_i^*\, ds_j) &=&
-\, s_i^*\,h\, s_j\,+\, s_i^*\,h\, s_j\ =\ 0\ ,\cr
\theta_h\lrcorner \sum_{i=1,\dots,n} (ds_i\, s_i^*\,+\, s_i\, ds_i^*) &=&
\sum_{i=1,\dots,n} (h\,s_i\, s_i^*\,-\, s_i\, s_i^*\, h)\ =\ h\,-\, h\ =\ 0\ ,\cr
{\cal L}_{\theta_h}(ds_i^*\, s_j\,+\, s_i^*\, ds_j) &=&
ds_i^*\,h\, s_j\,-\, d(s_i^*\,h)\, s_j\,-\, s_i^*\, h\, ds_j\,+\, s_i^*\, d(h\, s_j)  \cr
&=& -\, s_i^*\, dh\, s_j\,+\, s_i^*\, dh\, s_j\ =\ 0\ ,\cr
{\cal L}_{\theta_h}\sum_{i=1,\dots,n} (ds_i\, s_i^*\,+\, s_i\, ds_i^*) &=&
\sum_{i=1,\dots,n} (d(h\,s_i)\, s_i^*\,-\, ds_i\, s_i^*\, h\,-\, s_i\, d(s_i^*\, h)\,+\, h\,s_i\, ds_i^*)
\cr &=& \sum_{i=1,\dots,n} ( \, dh\, s_i\, s_i^*\,-\, s_i\, s_i^*\,dh\,)\ =\ 0\ .
\end{eqnarray*} 
Now we choose $\omega=\sum_i ds_i\,ds_i^*$, and note that $d\omega=0$. If we choose
$h=s_k\, s_l^*$, then
\begin{eqnarray*}
\theta_{s_k\, s_l^*}\lrcorner \,\omega&=& \sum_i (s_k\, s_l^*\, s_i\,ds_i^*\,+\, ds_i\, s_i^*\, s_k\, s_l^*)\cr
&=& s_k\, ds_l^*\,+\, ds_k\, s_l^*\ =\ d(s_k\, s_l^*)\ .
\end{eqnarray*} 
The set of derivations spanned by $\theta_{s_k\, s_l^*}$ for $1\le k,l\le n$ is closed under commutator,
and we call it $V$. We see that the Hamiltonian element corresponding to the derivation
$\theta_{s_k\, s_l^*}$ is $s_k\, s_l^*$, and that the Poisson brackets are given by
\[
\{s_k\, s_l^*,s_r\, s_m^*\}\ =\ \theta_{s_k\, s_l^*}(s_r\, s_m^*)\ =\ 
\delta_{l,r}\, s_k\, s_m^*\,-\, \delta_{m,k}\, s_r\, s_l^*\ .
\]

\section{An example of tensor products and interactions} 
Consider the algebra $A=C^\infty(\Bbb R^2,M_2(\Bbb R))$, where we use coordinates $x,y$ for $\Bbb R^2$.
The calculus we use will be the standard tensor product one, i.e.\
\[
\Omega^n(\, C^\infty(\Bbb R^2)\tens M_2(\Bbb R)\,)\ =\ \bigoplus_{p+q=n}
\Omega^p(\, C^\infty(\Bbb R^2)\,)  \tens
\Omega^q(\,  M_2(\Bbb R)\,)\ ,
\]
with $d$ operator and multiplication given by
\begin{eqnarray}
d(\tau\tens\eta) &=& d\tau\tens\eta\,+\, (-1)^{|\tau|}\, \tau\tens d\eta\ ,\cr
(\tau\tens\eta)\,  (\tau'\tens\eta') &=& (-1)^{|\eta|\, |\tau'|}\, \tau\tau'\tens\eta\eta'\ .
\end{eqnarray} 
The $d$ operator on $M_2(\Bbb R)$ is the one we defined earlier (we use $d_{M_2}$ to avoid confusion later),
and the $d$ operator on $C^\infty(\Bbb R^2)$ is the usual one:
\[
d\tau\ =\ dx\,\frac{\partial \tau}{\partial x}\,+\, dy\,\frac{\partial \tau}{\partial y}\ .
\]
There are derivations on the algebra given, for $f:\Bbb R^2\to M_2(\Bbb R)$, by
\[
\theta(f)\ =\ \theta_x\, \frac{\partial f}{\partial x}\,+\, 
\theta_y\, \frac{\partial f}{\partial y}\,+\, [\theta_S,f]\ ,
\]
where $\theta_x$ and $\theta_y$ are real valued functions
times the identity matrix on $\Bbb R^2$, and $\theta_S$ is
an antisymmetric matrix valued function on $\Bbb R^2$. This has evaluations on the 1-forms given by
$\theta\lrcorner dx=\theta_x$, $\theta\lrcorner dy=\theta_y$ and $\theta\lrcorner dE_{ij}=[\theta_S,E_{ij}]$. 
 We shall take the 2-form
\[
\omega\ =\ dx\, dy\,+\,\frac12 \sum_{ij} dE_{ij}\, dE_{ij}\,+\, dx\,(dE_{12}-dE_{21})\ ,
\]
where we have added the last term to ensure some interaction between the vector field
and the antisymmetric matrix parts of the derivations. Then we calculate
\[
\theta\lrcorner \,\omega\ =\ \theta_x\, dy\,-\, \theta_y\, dx\,+\, d_{M_2} \theta_S\,+\,
\theta_x\, (dE_{12}-dE_{21})\,-\, dx\, [\theta_S,E_{12}-E_{21}]\ .
\]
The last term here vanishes, since in $M_2(\Bbb R)$ any antisymmetric matrix is a multiple
of $E_{12}-E_{21}$, so
\[
\theta\lrcorner \,\omega\ =\ \theta_x\, dy\,-\, \theta_y\, dx\,+\, d_{M_2} (\theta_S\,+\,
\theta_x\, (E_{12}-E_{21}))\ .
\]
It is now reasonably simple to see that $\omega$ is non-degenerate for all
derivations of the form we are considering. 

Given an element of the algebra $a\in C^\infty(\Bbb R^2,M_2(\Bbb R))$
we have
\[
da\ =\ \frac{\partial a}{\partial x}\, dx\,+\, 
\frac{\partial a}{\partial y}\, dy\,+\, d_{M_2}a\ ,
\]
so if $\theta\lrcorner \,\omega=da$ then $d_{M_2}(\theta_S\,+\,\theta_x\, (E_{12}-E_{21}))=d_{M_2}a$,
$\theta_x=\frac{\partial a}{\partial y}$ and $\theta_y=-\frac{\partial a}{\partial x}$. We see that
if we put $a(x,y)\,=\, T\, +\, f(x,y)\, I_2$, where $T$ is a constant antisymmetric matrix
and $f(x,y)$ is a real valued function, then $(X_a)_x\,=\, \frac{\partial f}{\partial y}\, I_2$,
 $(X_a)_y\,=\,-\,\frac{\partial f}{\partial x}\, I_2$
and $(X_a)_S\,=\, T\,-\, \frac{\partial f}{\partial y}\, (E_{12}-E_{21})$. Now we can calculate the Poisson
bracket
of two such Hamiltonian functions:
\begin{eqnarray*}
\{ T\, +\, f(x,y)\, I_2, R\, +\, g(x,y)\, I_2\} &=& 
\frac{\partial f}{\partial y}\,\frac{\partial g}{\partial x}\, I_2\,-\,
\frac{\partial f}{\partial x}\,\frac{\partial g}{\partial y}\, I_2\,+\,
[T\,-\, \frac{\partial f}{\partial y}\,(E_{12}-E_{21})\,,\, R\,+\, g\, I_2]  \cr
&=& \Big( \frac{\partial f}{\partial y}\,\frac{\partial g}{\partial x}\,-\,
\frac{\partial f}{\partial x}\,\frac{\partial g}{\partial y}  \Big)\, I_2\ .
\end{eqnarray*}

\end{document}